\documentclass[12pts]{amsart}
\usepackage{latexsym,amsthm,amssymb,amsmath,amsfonts,euscript,amstext}
\usepackage{graphicx,epsfig,epstopdf}
\newtheorem{theorem}{Theorem}[section]

\theoremstyle{definition}
\newtheorem{defi}[theorem]{Definition}

\newtheorem{prop}[theorem]{Proposition}

\theoremstyle{remark}
\newtheorem{remark}[theorem]{Remark}

\usepackage{color}

\usepackage[utf8]{inputenc}

\numberwithin{equation}{section}

\date{}

\newcommand{\R}{\mathbb{R}}

\begin{document}
\title{Well-Posedness and Stability of the Stochastic OGTT Model}
\author{Paul Bekima}

\address{ Paul Bekima\newline
Department of Mathematics\\
Morgan State University\\
Baltimore, MD 21251, USA}
\email{pabek1@morgan.edu}

\maketitle

\section*{Introduction}

\vspace{8 pt}

The mathematical modeling of the Glucose - Insulin interaction was kicked off by Bergman et al. with the minimum model ([14.], [15.]). One of the main objective of the mathematical studies in the field of diabetes dynamics is to provide efficient models amenable to the estimation of the parameters, eventually with the greatest accuracy.
According to Boutayeb [3.], mathematical models in diabetes are generally of two types: simple and comprehensive.\

. Simple: with a minimum of parameters. \

. Comprehensive: with a large number of parameters.\



The current paper falls under the scoop of comprehensive models. Indeed, it is based on the paper by Ha and al.\textit{ Estimating insulin sensitivity and $ \beta $-cell function from the oral glucose tolerance test: validation of a new insulin sensitivity and secretion (ISS) model, J. American Physiological Society,  (2024).} 
To the proposed model, we added a perturbation driven by a family of independent white noises.  \\

The paper is organized as follow: we first establish the existence and uniqueness of a global positive solution to the proposed stochastic model. Then follows the study the statistical stability of the model by examining the existence of an invariant measure to the perturbed system. Once done with the well-posedness and the stability of the stochastic model, we then examine a Maximum Likelihood Estimation (MLE) scheme to estimate the parameters involved in the model. The paper ends with a brief discussion on potential future developments follows by an appendix whose goal is to make the paper as self contained as possible.  . \\

 Throughout this paper, unless otherwise specified, $ \lbrace \Omega, \mathcal{F},(\mathcal{F}_{t})_{t \geqslant 0}, P \rbrace $ is a complete probability space with a filtration $ (\mathcal{F}_{t})_{t \geqslant 0} $ satisfying the usual conditions; we will usually say $ \Omega $ as a shorter denomination for  $ \lbrace \Omega, \mathcal{F},(\mathcal{F}_{t})_{t \geqslant 0}, P \rbrace $. $ B_{i}(t), t\geqslant 0 $, i = 1,...,5 denote seven scalar brownian motions defined on the complete probability space $ \Omega $. \\

 \newpage

\section{The Model}




\vspace{8 pt}

The model is given as follows \\

\begin{equation}
    \begin{cases}
      dG(t) = (OGTT(t) + HGP(S_{I},I) -( E_{G0} + mS_{I}I(t))G(t)) dt + \alpha_{1}G(t)dB_{1}(t) \\\

   dI(t) = ( \frac{\beta}{V}ISG(\sigma,G) - kI(t))dt + \alpha_{2}I(t)dB_{2}(t)\\
      
   d \beta(t) = \frac{(P(ISR) - A(M)) \beta}{\tau_{\beta}} dt + \alpha_{3}I(t)dB_{3}(t) \\
   
  d \gamma(t) = \frac{ \gamma_{\infty}(G) - \gamma}{\tau_{\gamma}} dt + \alpha_{4}I(t)dB_{4}(t) \\

  d \sigma(t) = \frac{ \sigma_{\infty}(ISR,M) - \sigma}{\tau_{\sigma}} dt + \alpha_{5}I(t)dB_{5}(t) \\

    \end{cases}\,.
 \end{equation} \\

 $ t \in [0, T] $. \\

 With initial conditions \begin{equation}
 (G(0), I(0),\beta(0), \gamma(0), \sigma(0)) \in \R_{+}^{5} . \\\
\end{equation} \\

We can rewrite (1.1), (1.2) as \\

\begin{equation}
dY(t) = f(t,Y(t))dt + g(Y(t))dB(t), \ Y(0) \in \R_{+}^{5}
\end{equation} \\

Where $ Y(t) =  (G(t), I(t),\beta(t), \gamma(t), \sigma(t))^{T} $. \\\

\subsection*{Description of the Coefficients}

OGTT(t) is the glucose rate of appearance ($ R_{a} $) during the OGTT, modeled by a piece-wise linear function. \\

$ OGTT = \frac{OGTT_{0}}{V_{G}} $, where $ V_{G} = BW \times \bar{V} $ is the volume of distribution of glucose, BW is the body weight, $ \bar{V} = 1.569 \ dL/kg $, and \\

$ OGTT_{0} = a_{i-1} + \frac{a_{i} - a_{i-1}}{t_{i} - t_{i-1}} (t - t_{i-1}), \ t_{i-1} < t < t_{i}, \ i = 1,2,3 $, \

$ OGTT_{0} = 0 $ elsewhere. \\

$ HGP(S_{I},I) $ represents hepatic glucose production, which depends on peripheral insulin sensitivity ($ S_{I} $). The model HGP is a decreasing function of I. \\

$ ISR(\sigma, G) $ is the insulin secretion rate per unit mass. \

V is the volume of distribution. \

\begin{equation}
ISR = \sigma \frac{(M + \gamma)^{kISR}}{\alpha_{ISR}^{kISR} + (M + \gamma)^{kISR}}
\end{equation}

Where M was assumed to be a sigmoidally increasing function of G. \\

\begin{equation}
M = \frac{G^{kM}}{\alpha_{M}^{kM} + G^{kM}}
\end{equation} \\

$ \gamma $ represents the effect of ATP-sensitive $ K^{+} \ \ K_{ATP}$ channel density to shift the glucose dependence of secretion. \\

$ \gamma_{\infty} $ is an increasing sigmoidal function of G, $ \tau_{\gamma} $ is the time constant. \\

The parameters P and A are chosen such that modest increases in G result in a net increase in $ \beta $ but large increases in G result in a net decrease in $ \beta $. \\

See [1.] Supplemental Material, and [2.] for a complete description of the parameters. \\

\begin{remark}

The system (1.1) is non-autonomous due to the presence of OGTT(t). Through a Girsanov change of probability, the system can be transformed to an autonomous one.\\

Indeed, let $ dC(t) = (\frac{OGTT(t)}{\alpha_{1} G(t)} dB_{1}(t),0,0,0,0)^{T} $ \\

Taking the Poisson's bracket, \\

$ \langle dC(t), dY(t) \rangle \ = \ \frac{OGTT(t)}{\alpha_{1} G(t)} dt $ \\

Taking $ D(t) = exp(\int_{0}^{t} \frac{OGTT(s)}{\alpha_{1} G(s)} dB_{1}(s) - \frac{1}{2} \int_{0}^{t} (\frac{OGTT(s)}{\alpha_{1} G(s)})^{2}ds ) $ \\

The new probability Q is defined as $ dQ = D(t)dP\vert_{\mathcal{F}_{t}} $ \\

Under Q, the system becomes

\begin{equation}
    \begin{cases}
      dG(t) =  HGP(S_{I},I) -( E_{G0} + mS_{I}I(t))G(t)) dt + \alpha_{1}G(t)d \tilde{B}_{1}(t) \\\

   dI(t) = ( \frac{\beta}{V}ISG(\sigma,G) - kI(t))dt + \alpha_{2}I(t)dB_{2}(t)\\
      
   d \beta(t) = \frac{(P(ISR) - A(M)) \beta}{\tau_{\beta}} dt + \alpha_{3}I(t)dB_{3}(t) \\
   
  d \gamma(t) = \frac{ \gamma_{\infty}(G) - \gamma}{\tau_{\gamma}} dt + \alpha_{4}I(t)dB_{4}(t) \\

  d \sigma(t) = \frac{ \sigma_{\infty}(ISR,M) - \sigma}{\tau_{\sigma}} dt + \alpha_{5}I(t)dB_{5}(t) \\

    \end{cases}\,.
 \end{equation} \\

 $ t \in [0, T] $. \\

 With initial conditions \begin{equation}
 (G(0), I(0),\beta(0), \gamma(0), \sigma(0)) \in \R_{+}^{5} . \\\
\end{equation} \\

Where $ d \tilde{B}_{1}(t) = dB_{1}(t) + \frac{OGTT(t)}{\alpha_{1} G(t)}dt $ \\

Clearly $ \langle d \tilde{B}_{1}(t), d \tilde{B}_{1}(t) \rangle \ = \ \langle dB_{1}(t), dB_{1}(t) \rangle \ = \ dt $

\end{remark}

\newpage

\section{Positive Global Solution}

\vspace{8 pt} 

\begin{prop} 

For any given initial value $ Y(0) \in \R_{+}^{5} $  and $ t \geqslant 0 $ The equation ((1.3) has a unique global positive solution $ Y(t) \in \R_{+}^{5} $ almost surely, i.e. the solution will remain in $ \R_{+}^{5} $ with probability 1. \\\
\end{prop}

\begin{remark}
All the equations of (1.1) are of the form \\

\begin{equation}
dZ(t) = [b(t)+ \Phi(Y(t))Z(t)]dt + \alpha Z(t)dB(t)
\end{equation}

and can thus be written 

\begin{equation}
Z(t) = exp(\Phi(Y(t))t+ \alpha B(t) - \frac{1}{2}t)[b(0) + \int_{0}^{t} exp(-\Phi(Y(s))s- \alpha B(s) + \frac{1}{2}s) b(s)ds] \\
\end{equation}
\end{remark} 

\begin{remark}[Mao, page 58] 

If the assumptions of the existence-and-uniqueness theorem hold on every finite interval [0,T] of $ [0,\infty) $, then equation (1.3) has a unique solution Y(t) on the entire interval $ [0,\infty) $. Such a solution is called a global solution.\\\
\end{remark}

Let's recall the classical conditions for the existence and the unicity of a Global Solution. \

Given the Stochastic Differential Equation (SDE)\\

 $ dX(t) = f(t,X(t))dt + g(t,X(t))dB(t) $ \\

Assume that there exists two positive constants $ \bar{K} $ and K such that: \\

\begin{itemize}

\item[(i)] (Lipschitz Condition) for all x, y $ \in \ \R^{d} $ and $ t \in [t_{0}, T] $\

\begin{equation}
\vert f(t,x)- f(t,y) \vert^{2} \vee \vert g(t,x)- g(t,y) \vert^{2} \leqslant \bar{K} \vert x - y \vert^{2} \\
\end{equation}

\item[(ii)]( Linear Growth Condition )for all (t,x) $ \in \ [t_{0}, T]\times \R^{d} $ and $ t \in [t_{0}, T] $\

\begin{equation}
\vert f(t,x) \vert^{2} \vee \vert g(t,x) \vert^{2} \leqslant K( 1 + \vert x  \vert^{2}) \\
\end{equation}\\

\end{itemize}

Then, there exists a unique solution Y(t) to the equation (1.3) in $ \mathcal{M}^{2}([t_{0}, T]; \R^{d}) $ \\

Where $ \mathcal{M}^{2}([t_{0}, T]; \R^{d}) $ is the family of processes $ f(t)_{t_{0} \leqslant t \leqslant T} $ such that \

 $ E \int_{t_{0}}^{T} \vert f(t) \vert^{2}dt < \infty $ \\

Also, $ a \vee b $ is the maximum of a and b. \\

Let's also recall that the condition (2.4) can be weakened as follows. \\

\newpage

\begin{theorem}
Assume that the local Lipschitz condition holds, but the linear growth condition is replaced with the following monotone condition: \

There exists a positive constant K such that for all $ (t,y) \ \in \ [t_{0}, T] \times \R^{d} $ \

$ y^{T}f(t,y) + \frac{1}{2} \vert g(t,y) \vert^{2} \leqslant K(1 + \vert y\vert^{2} ) $. \\

Then, there exists a unique solution Y(t) to the equation (1.3) in $ \mathcal{M}^{2}([t_{0}, T]; \R^{d}) $ \\\
\end{theorem}

\begin{proof} (of Proposition 2.1)\\

The local existence follows from the fact that all the coefficients are Lipschitz. \

From remark 2.2, the solution is positive provided that b(0) positive; which is the case in the current model. \\

Theorem 2.4 and remark 2.3 ensure the global nature of the solution. \\

Indeed,\

 for the glucose equation \\

$ G(t)[OGTT(t) + HGP(S_{I},I) -( E_{G0} + mS_{I}I(t))G(t)] + \frac{1}{2} \alpha_{1}^{2} G^{2}(t) \leqslant (OGTT(t) + HGP(S_{I},I))^{2} + G^{2}(t) -E_{G0} G^{2}(t)  + \frac{1}{2} \alpha_{1}^{2} G^{2}(t) $ \\

OGTT is bounded on [0,T]. \\

HGP also is bounded, hence, there is $ K_{G} $ such that $ G(t)[OGTT(t) + HGP(S_{I},I) -( E_{G0} + mS_{I}I(t))G(t)] + \frac{1}{2} \alpha_{1}^{2} G^{2}(t) \leqslant K_{G}(1 + G^{2}(t)) $ \\

For I(t), since $ \frac{ISG((\sigma), G)}{V} $ is bounded, there is a $ M_{I} $ such that\

 $  \frac{\beta}{V}ISG(\sigma,G) - kI(t) \leqslant M_{I} \beta -kI(t) $ \
 
 thus $ I(t)( \frac{\beta}{V}ISG(\sigma,G) - kI(t)) \leqslant I(t)(M_{I} \beta -kI(t)) \leqslant K_{I}(1 + \beta(t)^{2} + I(t)^{2}) $ \
 
  for a well chosen $ K_{I} $.\\

  For $ \beta(t) $, $ \frac{(P(ISR) - A(M)) }{\tau_{\beta}} $ is bounded, thus there is $ K_{\beta} $ such that \
  
   $ \beta(t) \frac{(P(ISR) - A(M)) \beta(t)}{\tau_{\beta}} \leqslant (1+ K_{\beta}) \beta^{2} $. \\

 For $ \gamma(t) $ and $ \sigma(t) $, both $ \frac{\gamma_{\infty}}{\tau_{\gamma}} $ and $ \frac{\sigma_{\infty}(ISR,M)}{\tau_{\sigma}} $ are bounded and can thus be treated as in the case of $ \beta (t) $. \\

\end{proof}

\newpage

\section{Stability }

\vspace{8 pt} 

In this section, we are interested in the global stability of the dynamical system, unlike the Lyapunov stability which study stability properties of a single trajectory. \\

Let $ P(t,y,s,\cdot) $ the transition probability associated with $ Y(t) $ i.e \

$ P(t,y,s,\Gamma) = P(Y(s) \in \Gamma \vert Y(t) = y) $, assuming $ s > t $, also $ P(t,y,t,\cdot) = \delta_{y} $ the Dirac delta measure at y. \\

\begin{defi}
A probability measure $ \mu $ is invariant for the process Y(t) if \\

\begin{equation}
\int_{\Omega} d \mu (z) P(t,y,z,\Gamma) = \mu (\Gamma),\ \forall t \geqslant 0,  \ y \in \Omega  \,where\ \Omega \in \mathcal{B}(\R^{d})
\end{equation} \\

Let $ \mathcal{L}^{*} (\mu) = \partial_{x_{i}} \partial_{x_{j}} (a^{ij} \mu) - \partial_{x_{i}}(b^{i} \mu)  $ \\

The summation is taken over all repeated indices. \\

Where b = f(t,Y(t)), and $ a = g(t,Y(t))g(t,Y(t))^{T} $ as defined in equation (1.3). \\

The transition probability $ \mu = P(t,y,s,\cdot) $ which is the probability distribution of Y(t) satisfies the following \textit{Fokker-Planck-Kolmogorov equation} \\

\begin{equation}
 \partial_{t} \mu =  \partial_{x_{i}} \partial_{x_{j}} (a^{ij} \mu) - \partial_{x_{i}}(b^{i} \mu) 
\end{equation} \\

Any invariant distribution if it exist satisfies the \textit{stationary Fokker-Planck-Kolmogorov equation} \\

\begin{equation}
\partial_{x_{i}} \partial_{x_{j}} (a^{ij} \mu) - \partial_{x_{i}}(b^{i} \mu) = 0
\end{equation} \\

\end{defi}

\begin{prop}
There is a unique invariant measure to the system (equation 1.3) \\

\begin{equation}
dY(t) = f(t,Y(t))dt + g(Y(t))dB(t), \ Y(0) \in \R_{+}^{5}
\end{equation} \\

Where $ Y(t) =  (G(t), I(t),\beta(t), \gamma(t), \sigma(t))^{T} $. \\

\end{prop}

Before proceeding to the proof of the proposition, let's fix the notations and the language. \

Adopting the notations from [10.], let $ U_{R} = \lbrace x \in \R^{d}: \vert x \vert \leqslant R $; $ W^{d+,1}_{loc} (\R^{d}) $ denotes the class of all functions on $ \R^{d} $ whose restrictions to the balls $ U_{R} $ belongs to $ W^{p,1} (U_{R}) $ for some order $ p = p_{R} > d $. Similarly, $ L^{d+}_{loc}(\R^{d}) $ denotes the class of measurable functions whose restrictions to $ U_{R} $ belong to $ L^{p}(U_{R}) $ with some  $ p = p_{R} > d $. \\

We will use the following theorem to prove the existence of a solution to the equation (3.3), more over, the theorem establishes the existence of a density with respect to Lebesgue measure. \\

\begin{theorem} (Bogachev[10.], page 67)\\
 Let $ d \geqslant 2 $. Assume that for every ball $ U \in \R^{d} $ there exist numbers $ m_{U} > 0 $ and $ M_{U} > 0 $ such that \\
 
\begin{equation}
m_{U} \cdot I \leqslant A(x) \leqslant M_{U} \cdot I, \ \forall x \in U
\end{equation} \\

\begin{itemize}
\item[(i)] Suppose that $ a^{ij} \in W^{d+,1} (\R^{d}) $ and $ b^{i} \in L^{d+}_{loc} (\R^{d}) $. Then there exists a positive solution $ \varrho \in W^{d+,1} (\R^{d}) $ to the equation $ L^{*}_{A,b} (\varrho dx) =0 $. Moreover, if  $ a^{ij} \in W^{p,1} (\R^{d}) $ and $ b^{i} \in L^{p}_{loc} (\R^{d}) $ for some $ p > d $, then $ \varrho \in W^{p,1} (\R^{d}) $. \\

\item[(ii)] Suppose that $ a^{ij}, b^{i} \in L^{\infty}_{loc} (\R^{d}) $. Then there exists a nonnegative solution $ \varrho \in L^{d/(d-1)}_{loc} (\R^{d}) $ to the equation $ L^{*}_{A,b} (\varrho dx) =0 $ such that  \\

$ \int_{U_{1}} \varrho dx = 1  $. \\
\end{itemize}

Where $ A = (a^{ij})_{ij} $ and I is the identity matrix. \\
 
\end{theorem}

\begin{proof} We will proceed as follows \\

We apply the theorem to \\

\begin{equation}
dY^{\epsilon}(t) = f(t,Y^{\epsilon}(t))dt + [g(Y^{\epsilon}(t))+ \epsilon ]dB(t), \ Y^{\epsilon}(0) \in \R_{+}^{5}
\end{equation} \\ 

and then get $ \epsilon $ to tend to 0. \\

For any given $ \epsilon $, $ b = f(t,Y^{\epsilon}(t)) $ and $ a = g(Y^{\epsilon}(t))+ \epsilon $ meet the conditions for Theorem 3.3; hence there is and invariant measure $ \mu^{\epsilon} $ with a density $ \varrho^{\epsilon} $ such that \\

$ \int_{U_{1}} \varrho^{\epsilon} dx = 1  $. \\

Now, let's consider a sequence $ \epsilon_{n} $ that tends to 0 as n tends to $ \infty $, and let $ \varrho^{\epsilon}_{n} $ be the density associated with $ Y^{\epsilon_{n}} $.\\

We first show that, under the assumption that $ Y^{\epsilon_{n}}(0) \rightarrow Y(0) \ in \ L^{2}(P) $, we also have  $ Y^{\epsilon_{n}}(t) \rightarrow Y(t) \ in \ L^{2}(P) $ \\

$ E[\Vert Y^{\epsilon_{n}}(t) -  Y(t) \Vert^{2} ] \ \leqslant \  3E[\Vert Y^{\epsilon_{n}}(0) -  Y(0) \Vert^{2} ] \ + \ 3 E [ \Vert  \int_{0}^{t} (f(s,Y^{\epsilon_{n}}(s)) - f(s, Y(s)))ds \Vert^{2}] + \\
3   E [ \Vert  \int_{0}^{t} (g(s,Y^{\epsilon_{n}}(s))+ \epsilon - g(s, Y(s)))dB_{s} \Vert^{2}] $ \\

$ \leqslant  3E[\Vert Y^{\epsilon_{n}}(0) -  Y(0) \Vert^{2} ] \ + \ 3t E [   \int_{0}^{t} \Vert f(s,Y^{\epsilon_{n}}(s)) - f(s, Y(s))\Vert^{2} ds ] + \\
3 C  E [   \int_{0}^{t} \Vert g(s,Y^{\epsilon_{n}}(s))+ \epsilon - g(s, Y(s)) \Vert^{2} ds ] $ \\

The last inequality is obtained by using $ H \ddot{o}lder $ for the first integral and Burkholder-Davis-Gundy for the second one. \\

Using the Lipschitz condition, we obtain \\

$ E[\Vert Y^{\epsilon_{n}}(t) -  Y(t) \Vert^{2} ] \ \leqslant \  3E[\Vert Y^{\epsilon_{n}}(0) -  Y(0) \Vert^{2} + 6C \epsilon^{2}t] \ + \\ 3(\bar{K}t +  2\bar{C}) \int_{0}^{t} E[\Vert Y^{\epsilon_{n}}(t) -  Y(t) \Vert^{2} ]ds $ \\

Applying Gronwall, \\

$ E[\Vert Y^{\epsilon_{n}}(t) -  Y(t) \Vert^{2} ] \ \leqslant \  3E[\Vert Y^{\epsilon_{n}}(0) -  Y(0) \Vert^{2}] + 6C \epsilon^{2}t \ + \int_{0}^{t} 9E[\Vert Y^{\epsilon_{n}}(0) -  Y(0) \Vert^{2} + 6C \epsilon^{2}t] (\bar{K}t +  2\bar{C}) exp [\int_{s}^{t}3(\bar{K}t +  2\bar{C})du] ds $ \\

Hence if $ Y^{\epsilon_{n}}(0) \rightarrow Y(0) \ in \ L^{2}(P) $, then  $ Y^{\epsilon_{n}}(t) \rightarrow Y(t) \ in \ L^{2}(P) $ \\


Although this convergence solely dealt with the marginal distributions, given that we are involved with diffusion processes i.e. Markov processes with consistent finite dimensional distribution, the convergence is true for the whole processes. \\

Now we turn to the existence of an invariant measure with a density. \\

Let $ \mathcal{M}_{U_{1}}^{b} $ be the set of bounded positive measures on $ U_{1} $; \
equipped with the total variation norm (distance), $ \mathcal{M}_{U_{1}}^{b} $ is a topologically separated space.\
Since the total variation norm is lower semi-continuous, the sequence of measures $ \lbrace  \varrho^{\epsilon_{n}}dx \rbrace $ is relatively compact on $ \mathcal{M}_{U_{1}}^{b} $ thus tight since the space is separable and complete. (Let's recall from Theorem 3.3 that $ \int_{U_{1}} \varrho^{\epsilon_{n}} dx = 1 $  for all n).\
We can therefore extract a subsequence that converges to a measure $ \mu $ of the form $ \mu = \varrho dx $; since the associated cdfs $ F_{\epsilon_{n}}(t) $ are continuous and converge to the cdf $ F(t) $ of $ Y(t) $. \\

 It comes that, $ \mu = \varrho dx $ is an invariant measure to $ Y(t) $; (by Theorem 1.1, page 8 in Kifer [13])\\

 $ \mu $ is unique because $ Y(t) $ is Fellerian by Theorem 2.4. \\

\end{proof}

\newpage

\section{MLE estimations}

\vspace{8 pt}

Now that the existence, unicity and stability of the solutions to equations (1.3) has been established, it should be instructive to directly use the model to estimate the coefficients and compare the result to the deterministic model. This can be done using a maximum likelihood estimation (MLE) method. \

Here, we will approximate and discuss the existence and uniqueness of a consistent likelihood function for (1.3). \\



Let $ \theta $ be the collection of parameters in the system (1.3) (See the supplemental Material (Equations and Parameter Tables) for a full description of all the parameters).\

 Recall that $ Y(t) =  (G(t), I(t),\beta(t), \gamma(t), \sigma(t))^{T} $.\
 
  Let $ h(Y(t), \theta) = (f(t,Y(t), \theta) $, $B(t) = (B_{1}(t),...,B_{5}(t))^{T} $,\
  
   and $ \phi(Y(t),\theta)= g(Y(t), \theta) $  \\

Then the system (1.3) can be represented as \\

$ dY(t) = h(Y(t),\theta)dt + \phi(Y(t),\theta)dB(t) $ \\

As $ \Delta t \rightarrow 0 $,the Euler scheme produces the following discretization: \\

$ Y(t + \Delta t) - Y(t) = h(Y(t),\theta) \Delta t + \phi (Y(t),\theta)(B(t + \Delta t) - B(t)) $ \\\

Since $ B_{i}, \ i = 1,...,5 $ are independent Brownian motions, $ Y(t + \Delta t) - Y(t) $ are a five-dimensional Gaussian independent variables with mean $ h(Y(t),\theta) \Delta t $ and covariance matrix \\\

$$  
   A(Y(t),\theta)  =    \begin{bmatrix} \alpha_{1}^{2}G^{2}(t) & 0 & 0 & 0& 0\\ 0 & \alpha_{2}^{2}I^{2}(t) & 0 & 0 & 0 \\  0 & 0 & \alpha_{3}^{2} \beta^{2}(t) & 0 & 0 \\ 0 & 0 & 0 & \alpha_{4}^{2} \gamma^{2}(t) & 0\\ 0 & 0 & 0 & 0& \alpha_{5}^{2} \sigma^{2}(t) \end{bmatrix}
    $$ \\

Let $ \delta_{i}(\theta) = Y(t_{i}) - Y(t_{i-1}) - h(Y(t_{i-1}),\theta)(t_{i} - t_{i-1})$, where $ Y(t_{0}) = y_{0} $. We obtain the following log-likelihood function: \\\

$ ln(\theta \vert Y(t_{1}),..., Y(t_{n})) = - \frac{1}{2} \sum_{i=2}^{n} (\delta_{i}(\theta)^{\intercal}(\Sigma (Y(t_{i-1},\theta)(t_{i} - t_{i-1}))^{-1} \delta_{i}(\theta) + ln(\vert \Sigma (Y(t_{i-1},\theta)(t_{i} - t_{i-1})) \vert)) $ \\

By maximizing the above log-likelihood function with respect to $ \theta $, we  obtain a consistent estimator, denoted by $ \hat{\theta}_{n} $. \\

Now we left to assess the asymptotic behavior of $ \hat{\theta}_{n} $.

\newpage

\begin{prop} 

Under the following assumptions \\

\begin{itemize}

\item[A1] The diffusion matrix $ A(x; \theta) $ is positive definite \\

\item[A2] For each integer $ k \geqslant 1 $, the $ k^{th} $ order derivatives in x of the functions $ h (x; \theta) $ and $ \phi (x; \theta) $ exists and are locally bounded. \\

\item[A3] The transition density $ P_{X} (\Delta, x \vert x_{0}; \theta) $ is continuous in $ \theta $ and the log-likelihood function admits a unique maximizer i the parameter set. \\\

\end{itemize}

 $ \hat{\theta}_{n} $ converges in $ L^{2}(P) $ to a consistent estimator $ \hat{\theta} $.\\\

\end{prop}

\begin{remark}
Assumption A2 could be replaced with locally bounded derivatives and linear growth condition. \\\
\end{remark}

\begin{proof} (of Proposition 4.1)\\

It is sufficient to prove that $ lim_{n \rightarrow \infty} \sum_{i=0}^{2^{n}-1} \phi (Y(t_{i}),\theta)(B(t_{i+1}) - B(t_{i})) $ \

 exists. \\

Where the $ t_{i}'s  $ are the points of a subdivision of [0,T], with $ t_{i} = \frac{i}{2^{n}}T $ \\

Indeed, since the $ ( B(t_{i+1}) - B(t_{i}))_{i} $ are independent,\\

$ E[\sum_{i=0}^{2^{n}-1} \phi (Y(t_{i}),\theta)(B(t_{i+1}) - B(t_{i}))]^{2} \ \leqslant \ [\sum_{i=0}^{2^{n}-1} E[ \phi (Y(t_{i}),\theta)^{2}] \frac{(i+1)T}{2^{n}} $ \\

Since $ \phi (Y(t_{i}),\theta) $ is bounded, the sequence $  \sum_{i=0}^{2^{n}-1} \phi (Y(t_{i}),\theta)(B(t_{i+1}) - B(t_{i}) $ converges in $ L^{2}(P) $.\\

On the other hand, by the definition of the stochastic integral, that limit is almost surely equal to $ \int_{0}^{T} \phi(Y(t), \theta )dB(t) $. \\\
\end{proof}

\newpage

\section{Discussion}

\vspace{8 pt}

Through the fitting of ordinary differential equations, Ha and al. introduced in [1.] a new method to estimate insulin sensitivity and beta-cell function, then established the pertinence of their results by comparing their estimates to the ones from the ISS model with commonly used algebraic indices for OGTTs, and also to other differential equation-based fittings methods. Comparing the estimates obtained from the stochastic model to the ones obtained by Ha and al. will be decisive to the pertinence of the stochastic model. \\

Although for the MLE schemes many approaches are already available, the Euler scheme, the Monte-Carlo scheme, the Yoshida approach to name some of them, the toolbox can be enriched with numeric method for Partial Differential Equations (PDEs) if one handles the stochastic model through a martingale approach. \\

 The existence of an in variant measure enables a better controle of the dynamic of the system and eventually could help for a better understanding of the relationships between coefficients. If the initial conditions are drawn according to the invariant distribution, in a finite time, the dynamic will evolve to a point where it will be safe to replace the transition probability by the invariant measure. It becomes relevant to compute the time this substitution will happen. \\

Another direction is to consider other models for the noise; the noise could take into account the interaction among the variables. One should also consider the use of fractional brownian motions.\\\

\newpage

\section*{Appendice}

\vspace{8 pt}

\subsection*{A.1 Feller Processes} (Revuz -Yor PP 88 - 90)\\

Let $ C_{0}(E) $ be the space of continuous function on E that we simply write $ C_{0} $, and $ P_{t} $ a transition function. \\

\begin{prop} A transition function is Feller if and only if \\

\begin{itemize}
\item[(i)] $ P_{t} (C_{0}) \subset C_{0} $ for each t; \\

\item[(ii)] $ \forall f \in\in C_{0}, \ \forall x \in E, \ lim_{t \searrow 0} P_{t} f(x) = f(x) $ \\
\end{itemize}
\end{prop}

\subsection*{A.2 The Cameron-Martin-Girsanov Formula} (Stroock-Varadhan page 154)\\

\begin{theorem}

Let a: $ [0, \infty) \times \R^{d} \rightarrow S_{d} $, b: $ [0, \infty) \times \R^{d} \rightarrow \R^{d} $  and c: $ [0, \infty) \times \R^{d} \rightarrow \R^{d} $ be measurable functions such that a,b, and $ \langle c, ac\rangle $ are bounded. Then for each $ (s,x) \in [0, \infty) \times \R^{d} $ there is a one to one correspondence between solutions P to the martingale problem for a and b starting from (s,x) and solutions Q for a and b + ac starting from (s,x). The correspondence is such that $ P \rightarrow Q $ where $ Q \ll P $ on $ \mathcal{M}_{t}, \ t \geqslant 0 $, and dQ/dP on $ \mathcal{M}_{t} $ equals \\

$ exp [ \int_{s}^{t \veebar s} \langle c(u,x(u)), d \bar{x} (u) \rangle - \frac{1}{2} \int_{s}^{t \veebar s} \langle c(u),x(u)), a(u, x(u)) c(u, x(u))    \rangle du $ \\

with $ \bar{x}(t) = x(t) - \int_{s}^{t \veebar s} b(u, x(u))du, \ t \geq 0 $ \\\

\end{theorem}

\newpage

\section*{References}

\vspace{8 pt}

\begin{enumerate}

\item[1.] Ha J., Chung S.T., and al. Estimating insulin sensitivity and $ \beta $-cell function from the oral glucose tolerance test: validation of a new insulin sensitivity and secretion (ISS) model, \textit{J. American Physiological Society },  (2024). \\

\item[2.] Ha J., Sherman A., Type 2 diabetes: one disease, many pathways, Am. J. Physiol. Endocrinol Metab 319, E410 - E426, 2020. \\

\item[3.] Boutayeb A, Chetouani A, A critical review of Mathematical models and data used in diabetology, Biomedical Engineering Online, 2006. \\

\item[4.] Shi, X., Zheng, Q., Li, J., Zhou, X., Analysis of a stochastic IVGTT glucose-insulin model with time delay. (2020). \textit{AIMS, Mathematical Biosciences and Engineering}, 2310-2340. \\

\item[5.] Feigen P.D., Maximum Likelihood Estimation for Stochastic Processes - A Matingale Approach, PhD Dissertation, Australia National University, Sep. 1975. \\

\item[6.] Li C., Maximum Likelihood Estimation for Diffusion Processes via Closed-Form Density Expansions, The Annals of Statistics, Vol 41, No.3, 1350 - 1380, 2013. \\

\item[7.] Stroock D., Varadhan S., Multidimentional Diffusion Processes, Reprint of the 1997 Edition, Springer, 2006. \\

\item[8.] Khasminskii R., Stochastic Stability of Differential Equations, 2nd Edition. Springer, 2011. \\

\item[9.] Mao X, Stochastic Differential Equations and Applications, 2nd Edition, Woodhead Publishing,2011.\\

\item[10.] Bogachev V., Krylov N., Fokker-Planck-Kolmogorov Equations, Mathematical Surveys and Monographs, volume 207, 2015, American Mathematical Society.\\

\item[11.] Le Gall J.- F., Mouvement Brownien, Martingales et Calcul Stochastique. Springer, 2013.\\

\item[12.]  Revuz D. and Yor M., Continuous Martingales and Brownian Motion, 3rd Edition. Springer, 1999.\\

\item[13.] Kifer Y., Random Perturbations of Dynamical Systems, $ Birkh \ddot{a}user $, 1988. \\

\item[14.] Billingsley P., Convergence of Probability Measures, John Wiley \& Sons, NY, 1968. \\

\item[15.] Bergman R.N., Ider Y.Z., Bowden C.R.,Quantitative estimation of insulin sensitivity,\textit{ Am. J. Physiol. Endocrinol. Metab.}, \textbf{236} (1979), E667-667\\

\item[16.] Toffolo G., Bergman R.N., Finewood D.T.,Cobelli C.R.,  Quantitative estimation of beta cell sensitivity to glucose in the intact organism: A minimal model of insulin Kinetics in the dog. \textit{Diabetes}, \textbf{29} (1980), 979-990.\\

\end{enumerate}

\end{document}